\documentclass[11pt]{article}
\usepackage{amsmath,amsfonts}
\usepackage{verbatim}
\usepackage{latexsym}
\usepackage{graphicx}
\usepackage{dsfont}
\usepackage{indentfirst}
\usepackage{xcolor}
\usepackage{booktabs}
\makeatletter
\@addtoreset{equation}{section}

\begin{document}

\newcommand{\E}{\mathbb{E}}
\newcommand{\PP}{\mathbb{P}}
\newcommand{\RR}{\mathbb{R}}
\newcommand{\LL}{\mathbb{L}}

\newtheorem{theorem}{Theorem}[section]
\newtheorem{lemma}[theorem]{Lemma}
\newtheorem{corollary}[theorem]{Corollary}
\newtheorem{definition}[theorem]{Definition}
\newtheorem{assumption}[theorem]{Assumption}
\newtheorem{example}[theorem]{Example}
\newtheorem{proposition}[theorem]{Proposition}
\newtheorem{remark}[theorem]{Remark}

\newcommand\tq{{\scriptstyle{3\over 4 }\scriptstyle}}
\newcommand\qua{{\scriptstyle{1\over 4 }\scriptstyle}}
\newcommand\hf{{\textstyle{1\over 2 }\displaystyle}}
\newcommand\hhf{{\scriptstyle{1\over 2 }\scriptstyle}}

\newcommand{\proof}{\noindent {\it Proof}. }
\newcommand{\eproof}{\hfill $\Box$} % {\indent\vrule height6pt width4pt depth1pt\hfil\par\medbreak}

\def\a{\alpha} \def\g{\gamma}
\def\e{\varepsilon} \def\z{\zeta} \def\y{\eta} \def\ot{\theta}
\def\vo{\vartheta} \def\k{\kappa} \def\lbd{\lambda} \def\m{\mu} \def\n{\nu}
\def\x{\xi}  \def\r{\rho} \def\s{\sigma}
\def\p{\phi} \def\f{\varphi}   \def\w{\omega}
\def\q{\surd} \def\i{\bot} \def\h{\forall} \def\j{\emptyset}

\def\be{\beta} \def\de{\delta} \def\up{\upsilon} \def\eq{\equiv}
\def\ve{\vee} \def\we{\wedge}

\def\F{{\cal F}}
\def\T{\tau} \def\G{\Gamma}  \def\D{\Delta} \def\O{\Theta} \def\L{\Lambda}
\def\X{\Xi} \def\S{\Sigma} \def\W{\Omega}
\def\M{\partial} \def\N{\nabla} \def\Ex{\exists} \def\K{\times}
\def\V{\bigvee} \def\U{\bigwedge}

\def\1{\oslash} \def\2{\oplus} \def\3{\otimes} \def\4{\ominus}
\def\5{\circ} \def\6{\odot} \def\7{\backslash} \def\8{\infty}
\def\9{\bigcap} \def\0{\bigcup} \def\+{\pm} \def\-{\mp}
\def\la{\langle} \def\ra{\rangle}

\def\tl{\tilde}
\def\trace{\hbox{\rm trace}}
\def\diag{\hbox{\rm diag}}
\def\for{\quad\hbox{for }}
\def\refer{\hangindent=0.3in\hangafter=1}

\newcommand\wD{\widehat{\D}}

\title{
\bf   Consistency of least squares estimation to the parameter for stochastic differential equations under distribution uncertainty 
%\thanks{Partially supported by ...}
 }

\author{
{\bf
Chen Fei${}^{1}$,
Weiyin Fei${}^{2}$\thanks{Corresponding author. E-mail:  wyfei@ahpu.edu.cn}
%Litan Yan${}^{1}$
 }
\\
${}^1$Glorious Sun School of Business and Management,\\
Donghua University, Shanghai, 200051, China.\\
${}^2$ School of Mathematics and Physics, \\
Anhui Polytechnic University, Wuhu, Anhui 241000, China.\\
}
\date{}

\maketitle

\begin{abstract}
      Under distribution uncertainty, on the basis of discrete data we investigate the consistency of the least squares estimator (LSE) of the parameter for the stochastic differential equation (SDE) where the noise are characterized by $G$-Brownian motion. In order to obtain our main result of consistency of parameter estimation, we provide some lemmas by the theory of stochastic calculus of sublinear expectation. The result shows that under some regularity conditions, the least squares estimator is strong consistent uniformly on the prior set. An illustrative example is discussed.

\medskip \noindent
{\small\bf Key words: } stochastic differential equation disturbed by $G$-Brownian motiion ($G$-SDE); sublinear expectation; least squares estimator; exponential martingale inequality for capacity; strong consistency.

\medskip \noindent{\small \bf MR Subject Classification:} 60H10.
\end{abstract}
\section{Introduction}

In the study of statistic inferences of the parameters for stochastic differential equations, it has often
been supposed that the process in question can be observed continuously in time.
However, this assumption is not always appropriate in
practice. Even if it were possible to observe the process continuously,
in order to compute an estimate of the parameter one has to approximate the stochastic and
ordinary integrals by finite Riemann-Stieltjes sums. Therefore, it would be
interesting to know whether the parameter can be reasonably estimated from the
discrete data. For this kind of parameter estimation of the linear stochastic
differential equation disturbed by Brownian motion under certain probability distribution, we can refer to,  e.g., Bishwal \cite{Bi}, Kutoyants \cite{Ku}, and Prakasa Rao \cite{Pra}.
Based on the work of Hu and Nualart \cite{HN},  Shen et al. \cite{SYY} further discuss the least squares estimation for Ornstein-Uhlenbeck processes disturbed by the weighted fractional Brownian motion.

On the other hand, due to the well-known distinction attributed to Knight \cite{Kn}, there are two kinds of
uncertainty. The first, called risk, corresponds to situations in which all events relevant to
decision making are associated with obvious probability assignments. The second, called (Knightian) uncertainty, or ambiguity (following Ellsberg \cite{El}), corresponds to situations in which some events do not have an obvious probability assignment. Moreover, since the unobserved factor can affect the response randomly, Brownian motion has a conditional distribution
on the factor. We call this distribution uncertainty. A relevant example is in Huber \cite{Hub}. Often one  ignores the impact of
unobserved or ignored factors on the response. Under the framework of Knightian
uncertainty (ambiguity), different observations may come from different distributions
randomly selected from a class of distributions and the related problem analysis is based
on this distribution uncertainty (see, e.g., Lin {\sl et al.} \cite{Lin,Lin2}).

In fact, we know that the most vital one of all the assumption conditions imposed to classical statistical or probability models may be the distribution certainty
in the sense that each relevant random variable has a certain probability distribution that may or may not be known.
The classical linear expectation are based on such a key assumption. However, the distribution certainty
is not always the case in practice, such as in risk measure and robust portfolio decision-making in finance (for related references
see, e.g., Chen and Epstein \cite{CE},  Epstein and Ji \cite{EJ2}, Fei \cite{FeiSM}, Fei and Fei \cite{FF,FF1} and the references therein). Under distribution-uncertainty, the resulting expectation usually is nonlinear. The earlier works on sublinear expectation may
ascend to Huber \cite{Hub} for robust statistics, where he considered the static robust statistic analysis under uncertain distributions.
In the recent decades, the theory and methodology of nonlinear expectation have been well developed and received
much attention in some applied fields such as finance risk measure and portfolio making-decision (see, e.g., \cite{Ar, DPR}). Peng \cite{P2010} put forward recently the theory of sublinear expectation, where  the sublinear expectation is dynamic coherent. Under distribution uncertainty, Lin et al. \cite{Lin,Lin2} investigated upper expectation parameter regression and $k$-sample upper expectation linear regression modeling through using the idea of the nonconcave penalized likelihood suggested by Fan and Peng \cite{FP}.  Sun and Ji \cite{SJ} discussed the least squares estimator of random variables under sublinear expectations. Contrary to the fast development of the nonlinear probability theory, little attention was paid to the related
statistical inferences to the parameters for stochastic differential equations under distribution uncertainty to the best of our knowledge.

In reality, due to some complex environments, the system with ambiguity can be  characterized by stochastic differential equation driven by $G$-Brownian motion ($G$-SDE). Recently,  the related investigation is given by Fei et al. \cite{FFY},  Gao \cite{G}, Hu {\sl et al.} \cite{HJPS}, Lin\cite{Lq2013},  Lin \cite{Lin2013}, and Luo and Wang \cite{LW2014}, etc. As the classical case, we are faced with
the problem of the parameter estimator of $G$-SDE.

 In this paper, we discuss the $G$-SDE, where we need to estimate a unknown parameter in the drift term. By using sublinear expectation theory, we obtain the consistency of the parameter robust least squares estimation on the basis of discrete data, which is the main contribution of our paper. To our best knowledge, this is first to attempt discussing the topic of the parameter estimation of $G$-SDE under distribution uncertainty.

The main contribution of this paper is to establish the strong consistency of the parameter estimation to $G$-SDE under distribution-uncertainty by using the theory of sublinear expectations. For this, we prove the exponential martingale inequality under sublinear expectations and several key lemmas. The developed notions and
methodologies in our paper are nonclassical and original, and the theoretical framework establishes the foundations for
general nonlinear expectation statistics of stochastic differential equation driven by $G$-Brownian motion.

 The arrangement of the paper is as follows. In Section 2, we give preliminaries on sublinear expectations and $G$-Bwownian motions, then the related several lemmas are also given.
In Section 3, the strong convergence of the parameter of the $G$-SDE is investigated. For the proof of theorem, the related lemmas are further provided in Section 4. In Section 5, we discuss an example. Finally, we conclude in Section 6.

\section{ Preliminaries on sublinear expectation}

 In this section, we first give the notion of sublinear expectation space $(\Omega, {\cal H}, \hat{\mathbb{E}})$, where $\Omega$ is a given state set and $\cal H$ a linear space of real valued functions defined on $\Omega$. The space $\cal H$ can be considered as the space of random variables. The following concepts stem from  Peng \cite{P2010}.

\begin{definition} \label{Def2.1}
 A sublinear expectation $\hat{\mathbb E}$ is a functional $\hat{\mathbb  E}$: ${\cal H}\rightarrow {\mathbb  R}$ satisfying

\noindent(i) Monotonicity: $\hat{\mathbb{E}}[X]\geq \hat{\mathbb{E}}[Y]$ if $X\geq Y$;

\noindent(ii) Constant preserving: $\hat{\mathbb{E}}[c]=c$;

\noindent (iii) Sub-additivity: For each $X,Y\in {\cal H}$, $\hat{\mathbb{E}}[X+Y]\leq \hat{\mathbb{E}}[X]+\hat{\mathbb{E}}[Y]$;

\noindent(iv) Positivity homogeneity: $\hat{\mathbb{E}}[\lambda X]=\lambda \hat{\mathbb{E}}[X]$ for $\lambda\geq 0$.
\end{definition}

\begin{definition} \label{Def2.2} Let $(\Omega,{\cal H}, \hat{\mathbb{E}})$ be a sublinear expectation space. $(X(t))_{t\geq0}$ is called a $d$-dimensional stochastic process if for each $t\geq0$, $X(t)$ is a $d$-dimensional random vector in $\cal H$.

A $d$-dimensional process $(B(t))_{t\geq0}$ on a sublinear expectation space $(\Omega,{\cal H}, \hat{\mathbb{E}})$ is called a $G$-Brownian motion if the following properties are satisfies:

\noindent(i) $B_0(\omega)=0$;

\noindent (ii) for each $t,s\geq 0$, the increment $B(t+s)-B(t)$ is ${\cal N}(\{0\}\times s\Sigma)$-distributed and is independent from $(B(t_1),B(t_2),\cdots, B(t_n))$, for each $n\in {\mathbb N}$ and $0\leq t_1\leq \cdots\leq t_n\leq t$, where $\Sigma$ is a bounded, convex and closed subset of $d\times d$ nonnegative definite matrix family, and characterizes the covariance uncertainty of $G$-Brownian motion.

\end{definition}

More details on the notions of $G$-expectation ${\hat{\mathbb E}}$ and $G$-Brownian motion on the sublinear expectation space $(\Omega,{\cal H}, {\hat{\mathbb{E}}})$ may be found in Peng \cite{P2010}. Let $({\cal F}_t)_{t\geq 0}$ be a filtration generated by $G$-Brownian motion $(B(t))_{t\geq0}$. Define lower expectation ${\cal E}[X]:=\hat{\mathbb E}[X]$ for some $X\in {\cal H}$.

We now give the definition of It\^o integral. For the simplicity of expression, in the rest of the paper, all the processes take in one-dimensional real space ${\mathbb R}$. We introduce It\^o integral with respect to one dimensional $G$-Brownian motion with $G(\alpha):=\frac{1}{2}\hat{\mathbb{E}}[\alpha B(1)^2]=\frac{1}{2}(\bar{\sigma}^2\alpha^+-\underline{\sigma}^2\alpha^-)$, where $\hat{\mathbb{E}}[B(1)^2]=\bar{\sigma}^2, {\cal E}[ B(1)^2]=\underline{\sigma}^2$,
$0<\underline{\sigma}\leq \bar{\sigma}<\infty$.

Let $p\geq1$ be fixed. We consider the following type of simple processes: for a given partition $\pi_T=(t_0,\cdots,t_N)$
of $[0, T]$, where $T$ can take $\infty$, we get
$$\eta_t(\omega)=\sum_{k=0}^{N-1}\xi_k(\omega)I_{[t_k,t_{k+1})}(t),$$
where $\xi_k\in L^p_G({\cal F}_{t_k}), k=0,1,\cdots,N-1$ are given. The collection of these processes is denoted by $M_G^{p,0}(0,T)$. We denote by $M_G^p(0,T)$ the completion of $M_G^{p,0}(0,T)$  with the norm
 $$\|\eta\|_{M_G^p(0,T)}:=\left\{\hat{\mathbb{E}}\int_0^T|\eta(t)|^pdt\right\}^{1/p}<\infty.$$

We know that the weakly compact family of probability measures $\cal P$ characterizes the degree of Knightian uncertainty (see, e.g., Peng \cite{P2010}). Especially, if $\cal P$ is singleton, i.e. $\{P\}$, then the model has no ambiguity. The related calculus reduces to a classical one. We now define $G$-upper capacity $\mathbb{V}(\cdot)$ and $G$-lower capacity $\mathcal{V}(\cdot)$ by
\begin{align*}
{\mathbb V}(A)=&\sup_{P\in {\cal P}}P(A),~\forall A\in {\cal B}(\Omega),\\
{\mathcal V}(A)=&\inf_{P\in {\cal P}}P(A),~\forall A\in {\cal B}(\Omega).
\end{align*}
Thus a property is called to hold quasi surely (q.s.) if there exists a polar set $D$ with ${\mathbb V}(D_0)=0$ such that it holds for each $\omega\in D_0^c$.
A property is called to hold $\cal P$-q.s. if it holds for each $P\in {\cal P}$.
Let the generalized nonlinear expectation space be $(\Omega, {\mathcal H}, \hat{\mathbb E}, {\mathbb V}, ({\mathcal F}_t)_{t\geq0})$. Denote by $C_{l,Lip}({\mathbb R}^n)$ the space of all bounded
and Lipschitz real functions on ${\mathbb R}^n$.

\begin{definition}\label{D2.3} {\rm (1)}.  A real-valued process ($X(t)$, $t\in [0,\infty)$) with $X(t)\in L_G^2({\cal F}_t)$ for each $t\in [0,\infty)$ is called ergodic under sublinear expectation if the following holds
$${\mathcal V}\Big\{{\mathcal E}[X_0]\leq\inf\lim\limits_{T\rightarrow\infty}\frac{1}{T}\int_0^TX(t)dt\leq\sup\lim\limits_{T\rightarrow\infty}\frac{1}{T}\int_0^TX(t)dt\leq\hat{\mathbb E}[X_0]\Big\}=1.$$\\

{\rm (2)}. A real-valued process ($X(t)$, $t\in [0,T]$) with $X(t)\in L_G^2({\cal F}_t)$ for each $t\in [0,T]$ is called stationary under sublinear expectation  if for any $t_0,t_1,\cdots, t_n\in[0,T]$ with $t_0+t_i\in[0,T], i=\overline{1,n}$ we have $\hat{\mathbb E}[\varphi(X(t_1), \cdots,X(t_n))]=\hat{\mathbb E}[\varphi(X(t_0+t_1), \cdots,X(t_0+t_n))]$ for each $\varphi\in C_{l,Lip}({\mathbb R}^n)$.
\end{definition}

Obviously, if a real valued process $(X(t), t\in [0,T])$ is stationary, then we have, for any $t\in [0,T]$,
$$\hat{\mathbb E}[\varphi(X(t))]=\hat{\mathbb E}[\varphi(X(0))], ~\varphi\in C_{l,Lip}({\mathbb R}).$$
Specifically, we have $\hat{\mathbb E}[|X(t)|^p]=\hat{\mathbb E}[|X(0)|^p], ~p\geq 0.$

In what follows, we give the definition of $G$-martingale and the related inequalities. On a filtered sublinear expectation space $(\Omega, {\cal H}, ({\cal F}_t)_{t\in [0,T]}, {\widehat{\mathbb{E}}})$, we give the concept of $G$-martingales as follows.
\begin{definition}  A real-valued process ($X(t)$, $t\in [0,T]$) with $X(t)\in L_G^1({\cal F}_t)$ for each $t\in [0,T]$ is a $G$-(super, sub) martingale with respect to $({\cal F}_t)$ if $\widehat{\mathbb{E}}[X(t)\mid{\cal F}_s](\leq,\geq)= X_s$, q.s. for each pair $s, t$ such that $s<t$.

\end{definition}

\begin{lemma}\label{BC}
    {\rm(Borel-Cantelli Lemma (see, e.g., Chen \cite{Chen})} Let $\{A_n,n\geq1\}$ be a sequence of events in $\cal F$ and $({\mathbb V},{\mathcal V})$ be a pair of capacities generated by sublinear expectation  $\hat{\mathbb E}$. If $\sum\limits_{n=1}^\infty {\mathbb V}(A_n)<\infty,$ then ${\mathbb V}\left(\bigcap\limits_{n=1}^\infty\bigcup\limits_{k=n}^\infty A_k\right)=0$. That is, $\lim\sup\limits_{n\rightarrow\infty}A_n$ is a polar set.

\end{lemma}

\begin{lemma}\label{ExI}  {\rm(Exponential martingale inequality for capacity $\mathbb V$)} Let $g=(g_1,\cdots,g_m)\in{\mathcal L}^2_G({\mathbb R}_+;{\mathbb R}^m):=\{f: \mbox{for the stochastic process}~ f:~{\mathbb R}_+\rightarrow{\mathbb R}^m, ~\mbox{we have} ~ \hat{\mathbb E}\left[\int_0^\infty|f(t)|^2dt\right]<\infty\}$, and let $T, \alpha,\beta$ be positive numbers. Then we have
   \begin{align*}
     {\mathbb V}\Big\{\sup\limits_{0\leq t\leq T}\Big[\int_0^tg(s)dB(s)-\frac{\alpha}{2}\int_0^t|g(s)|^2ds\Big]>\beta\Big\}\leq e^{-\alpha\beta/\bar\sigma^2}.
     \end{align*}

\end{lemma}
{\textbf{ Proof}}:
Since $dB(t)=\sigma_tdW^{(\sigma_\cdot)}$, where $dW^{(\sigma_\cdot)}$ is a standard Brownian motion for each linear probability measure $P^{(\sigma_\cdot)}\in {\mathcal P}$ with an ${\cal F}_t$-adapted process $\sigma_t^2\in[\underline{\sigma}^2,\bar\sigma^2]$. From the classical exponential inequality (see, e.g., Theorem 1-7.4 in Mao \cite{Msde}) we get
 \begin{align*}
    &P^{(\sigma_\cdot)}\Big\{\sup\limits_{0\leq t\leq T}\Big[\int_0^tg(s)dB(s)-\frac{\alpha}{2}\int_0^t|g(s)|^2ds\Big]>\beta\Big\}\\
    &\leq P^{(\sigma_\cdot)}\Big\{\sup\limits_{0\leq t\leq T}\Big[\int_0^tg(s)\sigma_s dW^{(\sigma_\cdot)}(s)-\frac{\alpha}{2\bar\sigma^2}\int_0^t\sigma_s^2|g(s)|^2ds\Big]>\beta\Big\}\\
    &\leq e^{-\alpha\beta/\bar\sigma^2},
     \end{align*}
which easily shows the claim by the definition of the upper capcity $\mathbb V$.
 $\Box$

We now give the following the Burkholder-Davis-Gundy inequality (see, e.g., Gao \cite[Theorems 2.1-2.2]{G}).

\begin{lemma}(Burkholder-Davis-Gundy inequality)\label{BDG} Let $p\geq2$ and $\zeta=\{\zeta(s),s\in[0,T]\}\in M_G^p(0,T)$. Then, for all $t\in[0,T]$ such that
  \begin{align*}
 & \hat{\mathbb E}\sup\limits_{s\leq u\leq t}\Big|\int_s^u\zeta(v)d<B>(v)\Big|^p\leq(t-s)^{p-1}C_1(p,\bar\sigma)\hat{\mathbb E}\int_s^t|\zeta(v)|^pdv,\\
 &\hat{\mathbb E}\sup\limits_{s\leq u\leq t}\Big|\int_s^u\zeta(v)dB(v)\Big|^p\leq C_2(p,\bar\sigma)\hat{\mathbb E}\Big(\int_s^t|\zeta(v)|^2dv\Big)^{p/2},\\
\end{align*}
where the positive constants $C_i(p,\bar\sigma), i=1,2$ depend on parameters $p, \bar\sigma$.
\end{lemma}

\section{Consistency of least squares estimator}

Now we consider the linear stochastic differential equation
\begin{align} \label{SDE}
dX(t)= &a(\theta^0 ,X(t))dt +b(X(t)) d<B>(t)+dB(t), \quad X(0)=X_0,
\end{align}
where $X_0$ is a random variable, and $(B(t))_{t\geq0}$ is the $G$-Brownian motion in $\mathbb R$ on the generalized nonlinear expectation space $(\Omega, {\mathcal H}, \hat{\mathbb E}, {\mathbb V}, ({\mathcal F}_t)_{t\geq0})$.
  The functions $a(\cdot),b(\cdot)$ are assumed to know and  such that the solution of  Eq. (\ref{SDE}) exists  in $C([0,T],{\mathbb R})$.
   For Eq. (\ref{SDE}), we provide the following conditions.

   \begin{assumption} \label{A2.1}

   {\rm(A1)}. The parameter space $\Theta$ is a compact in $\mathbb R$.\\
   {\rm(A2)}. The functions $a(\theta,x),b(x)$ are Lipschizian, i.e., there exist nonnegative continuous functions $K_1(\cdot), D(\cdot)$ and a positive constant $K_2$ such that

   \begin{align*}
      &|a(\theta,x)-a(\theta,y)|\leq K_1(\theta) |x-y|,\\
   &|\rho(t)(b(x)-b(y))|\leq K_2|x-y|,\\
    &|a(\theta,x)-a(\psi,x)|\leq D(x)|\theta-\psi|
\end{align*}
    for all $\theta,\psi\in{\Theta}, x,y\in {\mathbb R}$, where $\rho(\cdot): [0,T]\rightarrow [\underline{\sigma}^2,\bar{\sigma}^2]$ is a continuous function, and
     \begin{align*}
   \sup\limits_{\theta\in \Theta}K_1(\theta)=K_1<\infty,\quad \hat{\mathbb E}|D(X_0)|^k=C_k<\infty, ~k>8.
\end{align*}
 We also have the following growth condition
 $$|a(\theta,x)|^2\vee|b(x)|^2\leq L(1+|x|^2)$$
 for some positive constants $L$.\\
 {\rm (A3)}. The process $(X(t), t\in[0,T])$ is stationary and ergodic under sublinear expectation (see, Definition \ref{D2.3}), and
  $$\hat{\mathbb E}|X_0|^k<\infty$$
 for some $k>8.$\\
 {\rm(A4)}. For the true parameter $\theta^0$ in (\ref{SDE}), we have
 $$\hat{\mathbb E}|a(\theta,X_0)-a(\theta^0, X_0)|^2=0\quad \mbox{iff}\quad \theta=\theta^0.$$
   \end{assumption}

 Next, our problem is to estimate the parameter $\theta^0$ on the basis of discrete observations
 $$X(t_0),X(t_1), \cdots,X(t_n), 0=t_0<t_1<\cdots<t_n=T.$$
Let $\Delta_it=t_i-t_{i-1}$. By the idea of the classical least squares estimation, we consider the following error
\begin{align} \label{LS}
{\mathbb  S}_{n,T}(\theta):=\sum\limits_{i=1}^n\frac{|X(t_i)-X(t_{i-1})-a(\theta,X({t_{i-1}}))\Delta_it-b(X(t_{i-1}))\sigma^2_{t_{i-1}}\Delta_it|^2}{\sigma^2_{t_{i-1}}\Delta_it},
\end{align}
where the distribution uncertainty results in $G$-Brownian noise with zero mean and uncertain quadratic variance $d<B>(t):=\sigma_t^2dt, \sigma_t^2\in[\underline\sigma^2,\bar\sigma^2]~q.s.$ We know that the values of ${\mathbb S}_{n,T}(\theta)$ in (\ref{LS}) change with $\sigma_t$ which characterize the probability measure $P^{(\sigma_t)}\in {\cal P}$. Now, define the following {\sl least squares estimator} (LSE)
\begin{align} \label{LSEa}
\hat{\theta}_{n,T}:=\arg\min\limits_{\theta\in \Theta}{\mathbb  S}_{n,T}(\theta),
\end{align}
where $\Theta$ is the space of parameter of the stochastic system with distribution uncertainty.

We will assume that $\Delta_it=t_i-t_{i-1}=T/n$ and $T=\Delta n^{1/2}$ for some fixed real number $\Delta>0$. However, we maintain the use of $T$ for convenience of notation. Let us decompose the sum of ${\mathbb S}_{n,T}(\theta)$ into three components
\begin{align*}
{\mathbb  S}_{n,T}(\theta)={\mathbb  S}_{n,T}(\theta^0)-2\sum\limits_{i=1}^n\frac{v_i\phi(\theta,X(t_{i-1}))}{\sigma^2_{t_i-1}}+\varphi_n^2(\theta),
\end{align*}
where
\begin{align*}
&v_i=X(t_i)-X(t_{i-1})-a(\theta^0,X(t_{i-1}))\Delta_it-b(X(t_{i-1}))\sigma^2_{t_{i-1}}\Delta_it,\\
&\phi(\theta,x)=a(\theta,x)-a(\theta^0,x),\\
&\varphi_n^2(\theta)=\sum\limits_{i=1}^n|\phi(\theta,X(t_{i-1}))|^2\Delta_it/\sigma^2_{t_{i-1}}.
\end{align*}
Now we define the following function
\begin{align} \label{Qn}
{\mathbb Q}_{n,T}(\theta):=T^{-1}({\mathbb S}_{n,T}(\theta)-{\mathbb S}_{n,T}(\theta^0))=T^{-1}\varphi_n^2(\theta)-2T^{-1}\sum\limits_{i=1}^n\sigma^{-2}_{t_{i-1}}v_i\phi(\theta,X(t_{i-1}).
\end{align}

We give the main result of this paper. The proof of theorem is further provided by several key lemmas which are placed in Section 4.

\begin{theorem} \label{Th}
Let Assumption \ref{A2.1} hold and
 \begin{align*}
&\sup\limits_{\theta\in\Theta}\Big|T^{-1}\sum\limits_{i=1}^n\frac{\phi(\theta,X(t_{i-1})}{\sigma^{2}_{t_{i-1}}}\Delta_i B\Big|\rightarrow 0\quad \mbox{q.s.}~ \mbox{as}~n\rightarrow\infty.
\end{align*}
Then $\hat\theta_{n,T}\rightarrow\theta^0$~q.s. as $n\rightarrow\infty$. That is, the least squares estimator is consistent to the true parameter $\theta^0$ of $G$-SDE (\ref{SDE}) quasi surely.
\end{theorem}
{\textbf{ Proof}} From Eq. (\ref{SDE}), we can have the following decomposition
\begin{align}\label{Tv}
&T^{-1}\sum\limits_{i=1}^n\frac{v_i\phi(\theta,X(t_{i-1}))}{\sigma^2_{t_{i-1}}}\nonumber\\
&= T^{-1}\sum\limits_{i=1}^n[X(t_i)-X(t_{i-1})-a(\theta^0,X(t_{i-1}))\Delta_it-b(X(t_{i-1}))\sigma^2_{t_{i-1}}\Delta_it]\frac{\phi(\theta,X(t_{i-1}))}{\sigma^2_{t_{i-1}}}\nonumber\\
&=T^{-1}\sum\limits_{i=1}^n\int_{t_{i-1}}^{t_i}a(\theta^0,X(s)-a(\theta^0,X(t_{i-1}))]\frac{\phi(\theta,X(t_{i-1}))}{\sigma^2_{t_{i-1}}}ds\nonumber\\
&\quad +T^{-1}\sum\limits_{i=1}^n\frac{\phi(\theta,X(t_{i-1}))}{\sigma^2_{t_{i-1}}}\Delta_iB,
\end{align}
where $\Delta_iB=B(t_i)-B(t_{i-1})$. By Lemma \ref{Lb}-\ref{LH}, we have, quasi surely,
$$T^{-1}\sum\limits_{i=1}^nv_i\frac{\phi(\theta,X(t_{i-1}))}{\sigma^2_{t_{i-1}}}\rightarrow0~\mbox{uniformly in}~\theta\in\Theta.$$
In terms of Lemma \ref{Lp} and (\ref{Qn}), we have, uniformly in $\theta$,
\begin{align*}
&{\mathcal V}\Big\{{\mathbb Q}_l(\theta)\leq \inf\lim\limits_{n\rightarrow\infty}{\mathbb Q}_{n,T}(\theta)\leq\sup\lim\limits_{n\rightarrow\infty}{\mathbb Q}_{n,T}(\theta)\leq{\mathbb Q}_u(\theta)\Big\}=1
\end{align*}
 where, the limiting functions ${\mathbb Q}_l(\theta)=\bar\sigma^{-2}{\mathcal E}|a(\theta,X_0)-a(\theta^0,X_0)|^2$  and ${\mathbb Q}_u(\theta)=\underline\sigma^{-2}\hat{\mathbb E}|a(\theta,X_0)-a(\theta^0,X_0)|^2$ fulfill the conditions (C2) and (C3) in Lemma \ref{LL} by Assumption \ref{A2.1} (A3). Thus, the proof of theorem is complete. $\Box$

\section{Several lemmas}
In this section, we provide several lemmas for the proof of the main result in Section 3.
For our aim, we now derive an estimate of $\mathbb E|X_t-X_s|^{2q}$ as a proposition. Now we discuss the more general
It\^{o} stochastic differential equation
\begin{align}\label{fg}
&dX(t)=f(t,X(t))dt+g(t,X(t))d<B>(t)+h(t,X(t))dB(t),~X(0)=x_0,~t\in[0,T],
\end{align}
where $(B(t))_{t\geq0}$ is one dimensional $G$-Brownian motion. We assume that the functions $f,g$ satisfy the following Lipschitz and growth conditions:

There is a constant $\bar K>0$ such that \\
(i) For all $ t\in[0,T], x\in {\mathbb R} $\\
$$|f(t,x)-f(t,y)|\vee|g(t,x)-g(t,y)|\vee|h(t,x)-h(t,y)|\leq \bar K|x-y|;$$\\
(ii) For all  $t\in[0,T],x \in {\mathbb R}$\\
$$|f(t,x)|^2\vee|g(t,x)|^2\vee|h(t,x)|^2\leq \bar K^2(1+|x|^2).$$

\begin{lemma}\label{Ip} If the above conditions (i), (ii)  and  Assumption \ref{A2.1} are satisfied, then for all $q\geq1, s,t\in [0,T] $ with $|t-s|<1,$ we have
$$\hat{\mathbb E}|X(t)-X(s)|^{2q}\leq K_3(t-s)^{q},$$
where the constant $K_3$ depends only on $q, \bar K, \bar\sigma$.
\end{lemma}

{\noindent\textbf{ Proof}} The growth condition (ii) implies
\begin{align}\label{Grow}
|f(t,x)|^{2q}\vee|g(t,x)|^{2q}\vee|h(t,x)|^{2q}&\leq \bar K^{2q}(1+|x|^2)^q \nonumber  \\
&\leq 2^{q-1}\bar K^{2q}(1+|x|^{2q}).
\end{align}
 By (\ref{fg}) and the H\"{o}lder inequality, we get
\begin{align}\label{pp}
&\hat{\mathbb E}|X(t)-X(s)|^{2q}\leq 2^{2q-1}\hat{\mathbb E}\Big\{\Big|\int_s^{t}f(u,X(u))du\Big|^{2q}\nonumber\\
&\quad+\Big|\int_s^{t}g(u,X(u))d<B>(u)\Big|^{2q}+\Big|\int_s^{t}h(u,X(u))dB(u)\Big|^{2q}\Big\}\nonumber\\
&\leq 2^{2q-1}\Big\{(t-s)^{2q-1}\int_s^{t}\hat{\mathbb E}|f(u,X(u))|^{2q}du\nonumber\\
&\quad+\hat{\mathbb E}\Big|\int_s^{t}g(u,X(u))d<B>(u)\Big|^{2q}+\hat{\mathbb E}\Big|\int_s^{t}h(u,X(u))dB(u)\Big|^{2q}\Big\}.
\end{align}
Due to $|t-s|<1$ and Assumption \ref{A2.1}, together with (\ref{Grow})-(\ref{pp}) and  the Burkholder-Davis-Gundy inequalities in Lemma \ref{BDG}, we have
\begin{align*}
&\hat{\mathbb E}|X(t)-X(s)|^{2q}\leq 2^{2q-1}\Big\{(t-s)^{2q-1}\int_s^{t}2^{q-1}\bar K^{2q}(1+\hat{\mathbb E}|X(u)|^{2q})du \nonumber\\
&\quad +C_1(q,\bar\sigma)(t-s)^{2q-1}\int_s^{t}\hat{\mathbb E}|g(u,X(u))|^{2q}du+C_2(q,\bar\sigma)\hat{\mathbb E}\Big(\int_s^{t}|h(u,X(u))|^2du\Big)^q\Big\}\\
&\leq 2^{2q-1}\Big\{(t-s)^{2q}2^{q-1}\bar K^{2q}(1+\hat{\mathbb E}|X_0|^{2q})+C_1(q,\bar\sigma)(t-s)^{2q-1}\int_s^{t}2^{q-1}\bar K^{2q}(1+\hat{\mathbb E}|X(u)|^{2q})du\\
&\quad +C_2(q,\bar\sigma)(t-s)^{q-1}\int_s^{t}2^{q-1}\bar K^{2q}(1+\hat{\mathbb E}|X(u)|^{2q})du\Big\}\\
&\leq 2^{3q-2}\bar K^{2q}(1+\hat{\mathbb E}|X_0|^{2q})(1+C_1(q,\bar\sigma)+C_2(q,\bar\sigma))(t-s)^q\\
&:=K_3(q,\bar K,\bar\sigma)(t-s)^q,
\end{align*}
where we also use the H\"older inequality on sublinear expectation (see, e.g., Peng \cite[Proposition 1.16 on page 95]{P2010}). Thus the proof is complete. $\Box$

\begin{lemma}\label{Lb} Let Assumption \ref{A2.1} hold. Then we have
\begin{align*}
\sup\limits_{\theta\in\Theta}\Big|T^{-1}\sum\limits_{i=1}^n\int_{t_{i-1}}^{t_i}[a(\theta^0,X(s))-a(\theta^0,X(t_{i-1}))]\frac{\phi(\theta,X(t_{i-1}))}{\sigma^2_{t_{i-1}}}ds\Big|\rightarrow0\quad\mbox{q.s.}
~\mbox{as}~T\rightarrow\infty.
\end{align*}
\end{lemma}
{\textbf{ Proof}} For $k>0$, we get
\begin{align}\label{W1}
&\hat{\mathbb E}\Big\{\sup\limits_{\theta\in\Theta}\Big|T^{-1}\sum\limits_{i=1}^n\int_{t_{i-1}}^{t_i}[a(\theta^0,X(s))-a(\theta^0,X(t_{i-1}))]\frac{\phi(\theta,X(t_{i-1}))}{\sigma^2_{t_{i-1}}}ds\Big|^{2k}\Big\}\nonumber\\
&\quad=\hat{\mathbb E}\Big\{\sup\limits_{\theta\in\Theta}\Big|T^{-1}\int_0^TY_n(s)ds\Big|^{2k}\Big\},
\end{align}
where $Y_n(s)=[a(\theta^0,X(s))-a(\theta^0,X(t_{i-1}))]\phi(\theta,X(t_{i-1}))/\sigma^2_{t_{i-1}}$ for $t_{i-1}\leq s<t_i$. Through the Cauchy-Schwartz inequality, from Assumption \ref{A2.1} (A2), we get that (\ref{W1}) is bounded by
\begin{align*}
&\leq T^{-2k}\hat{\mathbb E}\Big(\sup\limits_{\theta\in\Theta}T^{2k-1}\int_0^T|Y_n(s)|^{2k}ds\Big)\nonumber\\
&\leq T^{-2k}\hat{\mathbb E}\Big(\sup\limits_{\theta\in\Theta}T^{2k-1}\sum\limits_{i=1}^n\sigma^{-2}_{t_{i-1}}\int_{t_{i-1}}^{t_i}|a(\theta^0,X(s))-a(\theta^0,X(t_{i-1}))|^{2k}|\phi(\theta,X(t_{i-1}))|^{2k}ds\Big)\nonumber\\
&\leq T^{-1}\underline{\sigma}^{-2}Q_k\sum\limits_{i=1}^n\int_{t_{i-1}}^{t_i}\hat{\mathbb E}(|a(\theta^0,X(s))-a(\theta^0,X(t_{i-1}))|^{2k}|D(X(t_{i-1}))|^{2k})ds
\end{align*}
where $Q_k=:\sup\limits_{\theta\in\Theta}|\theta-\theta^0|^{2k}<\infty$. From the H\"older inequality on sublinear expectation, we have
\begin{align*}
&\leq T^{-1}\underline{\sigma}^{-2}Q_k\sum\limits_{i=1}^n\int_{t_{i-1}}^{t_i}(\hat{\mathbb E}|a(\theta^0,X(s))-a(\theta^0,X(t_{i-1}))|^{4k})^{1/2}(\hat{\mathbb E}D(X(t_{i-1}))|^{4k})^{1/2}ds\nonumber\\
&\leq T^{-1}\underline{\sigma}^{-2}Q_k(K(\theta^0))^{2k}(\hat{\mathbb E}D(X_0)^{4k})^{1/2}\sum\limits_{i=1}^n\int_{t_{i-1}}^{t_i}(\hat{\mathbb E}|X(s)-X(t_{i-1})|^{4k})^{1/2}ds,
\end{align*}
by Assumption \ref{A2.1} (A2) and (A3). Due to Assumption \ref{A2.1}, we know the conditions on coefficients of equation (\ref{fg}) hold. Thus, from Lemma \ref{Ip}, we get that (\ref{W1}) is bounded by
\begin{align*}
&\leq T^{-1}\underline{\sigma}^{-2}Q_k(K(\theta^0))^{2k}(\hat{\mathbb E}D(X_0)^{4k})^{1/2}K_3^{1/2}\sum\limits_{i=1}^n\int_{t_{i-1}}^{t_i}(s-t_{i-1})^kds\\
&\leq \underline{\sigma}^{-2}Q_k(K(\theta^0))^{2k}(\hat{\mathbb E}D(X_0)^{4k}K_3)^{1/2}T^{-1}\sum\limits_{i=1}^n\frac{(\Delta_it)^{k+1}}{k+1}\\
&\leq \frac{\underline{\sigma}^{-2}Q_k(K(\theta^0))^{2k}}{k+1}(\hat{\mathbb E}D(X_0)^{4k}K_3)^{1/2}\Delta^kn^{-k/2}, ~k>2.
\end{align*}
By Markovian inequality for upper capacity ${\mathbb  V}$ (see, e.g., Peng \cite{P2010}), Assumption \ref{A2.1},  and the above inequality, we obtain that, for $k>2$, $\forall \varepsilon>0$,
\begin{align*}
\sum\limits_{n=1}^\infty {\mathbb V}\Big\{\sup\limits_{\theta\in\Theta}\Big|T^{-1}\sum\limits_{i=1}^n\int_{t_{i-1}}^{t_i}[a(\theta^0,X(s))-a(\theta^0,X(t_{i-1}))]\frac{\phi(\theta,X(_{i-1}))}{\sigma^2_{t_{i-1}}}ds\Big|>\varepsilon\Big\}<\infty,
\end{align*}
by which, together with the Borel-Camtelli Lemma on sublinear expectation (see Lemma \ref{BC}), we get the desired result. Thus the proof is complete. $\Box$

\begin{lemma}\label{LH} Assume that for each $x,\phi(\cdot,x)$ takes value in a Hilbert space $\mathbb H$ with the property that $\mathbb H$ is continuously imbeded into the space $C(\Theta)$ of continous functions on $\Theta$. Suppose that $\{e_k\}$ is a complete orthonormal basis for $\mathbb H$ and that
\begin{align}\label{Re}
\phi(\theta,x)=\sum\limits_{k}\lambda_k(x)e_{k}(\theta)
\end{align}
satisfying
\begin{align}\label{kc}
|\lambda_k(x)|\leq c_k|x|~\mbox{and}~ \sum\limits_{k}k^{1+\beta}c_k^4<\infty
\end{align}
for some positive number $\beta$. Then we have
$$\sup\limits_{\theta\in\Theta}\Big|T^{-1}\sum\limits_{i=1}^n\sigma^{-2}_{t_{i-1}}\phi(\theta,X(t_{i-1}))\Delta_iB\Big|\rightarrow0\quad\mbox{q.s.}~\mbox{as}~T,n\rightarrow\infty.$$
\end{lemma}
{\textbf{ Proof}} Since $\mathbb H\hookrightarrow C(\Theta),$ we need only to show that
$$\Big\|T^{-1}\sum\limits_{i=1}^n\sigma^{-2}_{t_{i-1}}\phi(\theta,X(t_{i-1}))\Delta_iB\Big\|_{\mathbb H}\rightarrow0\quad\mbox{q.s.~as}~T,n\rightarrow\infty.$$
In view of expression (\ref{Re}), we have
 \begin{align*}
 \Big\|T^{-1}\sum\limits_{i=1}^n\sigma^{-2}_{t_{i-1}}\phi(\theta,X(t_{i-1}))\Delta_iB\Big\|_{\mathbb H}^2= \sum\limits_k\Big(T^{-1}\sum\limits_{i=1}^n\sigma^{-2}_{t_{i-1}}\lambda_k(X(t_{i-1}))\Delta_iB\Big)^2\quad\mbox{q.s.}
 \end{align*}
 Set
 $$\Lambda_{kn}=\sum\limits_{k}\sigma^{-2}_{t_{i-1}}\lambda_k(X(t_{i-1}))\chi_{[t_{i-1},t_i)}(s),$$
 where $\chi_{[t_{i-1},t_i)}, i=1,\cdots,n$, are indicative functions. Then, we get
  \begin{align}\label{La}
 T^{-1}\sum\limits_{i=1}^n\sigma^{-2}_{t_{i-1}}\lambda_i(X(t_{i-1})\Delta_iB=T^{-1}\int_0^T\Lambda_{kn}(s)dB(s)
 \end{align}
 and the integral on the RHS of (\ref{La})
 is well difined as $\Lambda_{kn}$ are adapted. By the exponential martingale inequality of Lemma \ref{ExI}, we obtain
 \begin{align*}
     {\mathbb V}\Big\{\int_0^T\Lambda_{kn}(s)dB(s)-\frac{\alpha}{2}\int_0^T\Lambda^2_{kn}(s)ds>\gamma\Big\}\leq e^{-\alpha\gamma/\bar\sigma^2}
\end{align*}
for any $\alpha,\gamma>0$. Thus we have
\begin{align}\label{VI}
     {\mathbb V}\Big\{\frac{1}{T}\int_0^T\Lambda_{kn}(s)dB(s)>\frac{\gamma}{T}+\frac{\alpha}{2T}\int_0^T\Lambda^2_{kn}(s)ds\Big\}\leq e^{-\alpha\gamma/\bar\sigma^2},\nonumber\\
     {\mathbb V}\Big\{\Big|\frac{1}{T}\int_0^T\Lambda_{kn}(s)dB(s)\Big|>\frac{\gamma}{T}+\frac{\alpha}{2T}\sum\limits_{i=1}^n\sigma^{-4}_{t_{i-1}}\lambda_k^2(X(t_{i-1}))\Delta_it\Big\}\leq2e^{-\alpha\gamma/\bar\sigma^2}.
\end{align}
By condition (\ref{kc}), we have
\begin{align*}
     T^{-1}\sum\limits_{i=1}^n\sigma^{-4}_{t_{i-1}}\lambda_i^2(X(t_{i-1}))\Delta_it\leq c_k^2T^{-1}\sum\limits_{i=1}^n\sigma^{-2}_{t_{i-1}}X^2(t_{i-1})\Delta_it.
\end{align*}
From Assumption \ref{A2.1} (A3), we know that $(X(t), t\geq 0)$ is a stationary and ergodic process under sublinear expectation. Moreover we have
\begin{align*}
    &{\mathcal V}\Big\{\bar\sigma^{-4}{\mathcal E}[X_0^2]\leq T^{-1}\inf\lim\limits_{T\rightarrow\infty}\sum\limits_{i=1}^nX^2(t_{i-1})\Delta_it\\
    &\quad\leq T^{-1}\sup\lim\limits_{T\rightarrow\infty}\sum\limits_{i=1}^nX^2(t_{i-1})\Delta_it\leq \underline\sigma^{-4}{\mathbb E}[X_0^2]\Big\}=1.
\end{align*}
Thus we know that there exists a random variable $\zeta$ with ${\mathcal V}(\zeta<\infty)=1$ such that
\begin{align*}
  T^{-1}\sum\limits_{i=1}^n\sigma^{-4}_{t_{i-1}}X^2(t_{i-1})\Delta_it\leq\zeta\quad\mbox{q.s. for all}~ T>0, n=1,2,\cdots.
\end{align*}
 Let
$$\Psi_{kn}:=\frac{1}{t_n}\int_0^{t_n}\Lambda_{kn}(s)dB(s).$$
Instead of $T$ by $t_n$ in order to stress the dependence of $T$ on $n$. We select
$$\alpha=\frac{k^\mu}{t_n^\delta},\quad \gamma=\frac{t_n^\eta}{k^\nu}$$
where $\delta<\eta<1$ and $1/2<\nu<\mu:=(1+\beta)/2$. For convenience of presentation, we can take enough large integer numbers $n,k$ such that $\frac{1}{t_n^{\eta-\delta}}+\frac{1}{k^{\mu-\nu}}\leq 1$, i.e., $t_n^{\eta-\delta}k^{\mu-\nu}\geq k^{\mu-\nu}+t_n^{\eta-\delta}$. Hence, we have
\begin{align*}
{\mathbb V}\Big(|\Psi_{kn}|>\frac{1}{t_{n}^{1-\eta}k^{\nu}}+\frac{k^\mu c_k^2\zeta}{2t_n^\delta}\Big)\leq2\exp(-k^{\mu-\nu}t_n^{\eta-\delta}/\bar\sigma^2),
\end{align*}
which shows
\begin{align*}
&{\mathbb V}\Big(\sum\limits_{k}\Psi_{kn}^2>\sum\limits_{k}\Big(\frac{1}{t_n^{1-\eta}k^\nu}+\frac{k^\mu c_k^2\zeta}{2t_n^\delta}\Big)^2\Big)\\
&\leq\sum\limits_{k}{\mathbb V}\Big(\Psi_{kn}^2>\Big(\frac{1}{t_n^{1-\eta}k^\nu}+\frac{k^\mu c_k^2\zeta}{2t_n^\delta}\Big)^2\Big)\\
&=\sum\limits_{k}{\mathbb V}\Big(|\Psi_{kn}|>\frac{1}{t_n^{1-\eta}k^\nu}+\frac{k^\mu c_k^2\zeta}{2t_n^\delta}\Big)\\
&\leq2\sum\limits_k\exp(-k^{\mu-\nu}t_n^{\eta-\delta}/\bar\sigma^2)\\
&\leq2\exp(-t_n^{\eta-\delta}/\bar\sigma^2)\sum\limits_k\exp(-k^{\mu-\nu}/\bar\sigma^2).
\end{align*}
Moreover, we get
\begin{align*}
&\sum\limits_n{\mathbb V}\Big(\sum\limits_k\Psi_{kn}^2>\sum\limits_k\Big(\frac{1}{t_n^{1-\eta}k^\nu}+\frac{k^\mu c_k^2\zeta}{2t_n^\delta}\Big)^2\Big)\\
&\leq2\sum\limits_n\exp(-t_n^{\eta-\delta})\sum\limits_k\exp(-k^{\mu-\nu})<\infty,
\end{align*}
as $\eta-\delta>0$ and $\mu-\nu>0. $ From the above, we get
\begin{align*}
\sum\limits_n{\mathbb V}\Big(\sum\limits_k\Psi_{kn}^2>\frac{2}{t_n^{2(1-\eta)}}\sum\limits_kk^{-2\nu}+\frac{\zeta^2}{t_n^{2\delta}}\sum\limits_kk^{2\mu}c_k^4\Big)<\infty.
\end{align*}
Since, by condition (\ref{kc}), $\delta<\eta<1$ and $1/2<\nu<\mu:=(1+\beta)/2$, we easily deduce
\begin{align*}
\frac{2}{t_n^{2(1-\eta)}}\sum\limits_kk^{-2\nu}+\frac{\zeta^2}{t_n^{2\delta}}\sum\limits_kk^{2\mu}c_k^4\rightarrow0\quad{\rm q.s.}~\mbox{as}~n\rightarrow\infty.
\end{align*}
From Lemma \ref{BC} and the assumption of the theorem, we get
$$\sum\limits_k\Big(t_n^{-1}\sum\limits_{i=1}^n\sigma^{-2}_{t_{i-1}}\lambda_k(X(t_{i-1}))\Delta_iB\Big)^2\rightarrow0\quad\mbox{q.s. as}~n\rightarrow\infty.$$
Thus, the proof is complete. $\Box$

Now, we prove the following claim.

\begin{lemma}\label{Lp} Under the Assumptions \ref{A2.1}, we have
\begin{align*}
&{\mathcal V}\Big\{\bar\sigma^{-2}{\mathcal E}|\phi(\theta,X_0)|^2\leq \inf\lim\limits_{T\rightarrow\infty}T^{-1}\varphi_n^2(\theta)\nonumber\\
&\quad\leq\sup\lim\limits_{T\rightarrow\infty}T^{-1}\varphi_n^2(\theta)\leq\underline\sigma^{-2}\hat{\mathbb E}|\phi(\theta,X_0)|^2\Big\}=1
\end{align*}
\end{lemma}

uniformly in $ \theta$.

{\noindent\textbf{ Proof}} By Assumption \ref{A2.1} (A3), we deduce
\begin{align}\label{VV}
&{\mathcal V}\Big\{\bar\sigma^{-2}{\mathcal E}|\phi(\theta,X_0)|^2\leq \inf\lim\limits_{T\rightarrow\infty}T^{-1}\int_0^T\tilde{\sigma}_t|\phi(\theta,X_s)|^2ds\nonumber\\
&\leq\sup\lim\limits_{T\rightarrow\infty}T^{-1}\int_0^T\tilde{\sigma}_t|\phi(\theta,X_s)|^2ds\leq\underline\sigma^{-2}\hat{\mathbb E}|\phi(\theta,X_0)|^2\Big\}
\end{align}
for each $\theta\in\Theta$ and $\tilde{\sigma}_t=\sum_i\sigma_{t_{i-1}}^{-2}\chi_{[t_{i-1},t_i)}(t),t\in[0.T]$. From Assumption  \ref{A2.1} (A1)-(A3) we get
\begin{align*}
&T^{-1}\int_0^T|\phi(\theta,X_s)|^2ds \leq T^{-1}|\theta-\theta^0|^2\int_0^T|D(X_s)|^2ds \nonumber  \\
&\leq\sup_{\theta\in\Theta}|\theta-\theta^0|^2T^{-1}\int_0^T|D(X_s)|^2ds \nonumber \\
&\leq \xi_1,
\end{align*}
quasi surely for some finite positive random variable $\xi_1$. From Assumption  \ref{A2.1} we also deduce that
\begin{align*}
\Big|T^{-1}\int_0^T\tilde{\sigma}_t|\phi(\theta_1,X_s)|^2ds-T^{-1}\int_0^T\tilde{\sigma}_t|\phi(\theta_2,X_s)|^2ds\Big|\leq \xi_2|\theta_{1}-\theta_{2}|
\end{align*}
quasi surely for some finite random variable $\xi_2$ and $\theta_1,\theta_2\in\Theta$. Hence the family of functions
$$\Big\{T^{-1}\int_0^T\tilde{\sigma}_t|\phi(
\cdot,X_s)|^2ds,T\geq 0\Big\}$$
are equicontinuous, which shows, by the Arzel\`{a}-Ascoli theorem, the convergence (\ref{VV}) is
uniform. Next we show that
\begin{align*}
T^{-1}\int_0^T\tilde{\sigma}_t|\phi(\theta,X_s)|^2ds-T^{-1}\varphi_n^2(\theta)\rightarrow 0
\end{align*}
quasi surely uniformly in $\theta$. Hence  we get
\begin{align*}
&\hat{\mathbb E}\Big\{\sup_{\theta\in\Theta}\Big|\int_0^T\tilde{\sigma}_t|\phi(\theta,X_s)|^2ds-\varphi_n^2(\theta)\Big|^{2q}\Big\}\nonumber\\
&=\mathbb E\Big\{\sup_{\theta\in\Theta}\Big|\int_0^T\tilde{\sigma}_t|\phi(\theta,X_s)|^2ds-\sum_{i=1}^n\sigma^{-2}_{t_{i-1}}|\phi(\theta,X(t_{i-1}))|^2\Delta_{i}t\Big|^{2q}\Big\}  \nonumber  \\
&=\mathbb E\Big\{\sup_{\theta\in\Theta}\Big|\sum_{i=1}^n\int_{t_{i-1}}^{t_{i}}\sigma^{-2}_{t_{i-1}}[\phi(\theta,X_s)-\phi(\theta,X(t_{i-1}))][\phi(\theta,X_s)+\phi(\theta,X(t_{i-1}))]ds\Big|^{2q}\Big\}.
\end{align*}
Using the H\"{o}lder inequality on sublinear expectation, we have
\begin{align*}
&\leq \underline\sigma^{-4q}T^{2q-1}\hat{\mathbb E}\Big\{\sup_{\theta\in\Theta}\sum_{i=1}^n\int_{t_{i-1}}^{t_i}|\phi(\theta,X_s)-\phi(\theta,X(t_{i-1}))|^{2q}|\phi(\theta,X_s)\\
&\qquad+\phi(\theta,X(t_{i-1}))|^{2q}ds\Big\} \\
&\leq \underline\sigma^{-4q}T^{2q-1}\sum_{i=1}^n\int_{t_{i-1}}^{t_i}\hat{\mathbb E}[\sup_{\theta\in\Theta}|\phi(\theta,X_s)-\phi(\theta,X(t_{i-1}))|^{2q}\sup_{\theta\in\Theta}|\phi(\theta,X_s)\\
&\qquad+\phi(\theta,X(t_{i-1}))|^{2q}]ds\\
&\leq \underline\sigma^{-4q}T^{2q-1}K^{2q}2^{2q-1}D_q\sum_{i=1}^n\int_{t_{i-1}}^{t_i}\hat{\mathbb E}[|X_s-X(t_{i-1})|^{2q}(|D(X_s)|^{2q}\\
&\qquad+|D(X(t_{i-1}))|^{2q})]ds \\
&\leq \underline\sigma^{-4q}T^{2q-1}K^{2q}2^{2q}D_q\sum_{i=1}^n\int_{t_{i-1}}^{t_i}(\hat{\mathbb E}|X_{s}-X(t_{i-1})|^{4q})^{1/2}(\hat{\mathbb E}|D(X_s)|^{4q}\\
&\qquad+\hat{\mathbb E}|D(X(t_{i-1}))|^{4q})^{1/2}ds,
\end{align*}
where $D_q=\sup_{\theta\in\Theta} |\theta-\theta^0|^{2q}< \infty$.  By Lemma \ref{Ip}, we get
\begin{align*}
&\leq \underline\sigma^{-4q}T^{2q-1}K^{2q}2^{2q+1}D_qK_3^{1/2}(\mathbb E|D(X_0)|^{4q})^{1/2}\sum_{i=1}^n\int_{t_{i-1}}^{t_i}(s-t_{i-1})^{q}ds \\
&\leq L_qT^{2q-1}n(T/n)^{q+1}, \\
\end{align*}
where $L_{q}=\underline\sigma^{-4q}K^{2q}2^{2q+1}D_qK_3^{1/2}(\hat{\mathbb E}|D(X_0)|^{4q})^{1/2}$. Thus, for $q>2$,
\begin{align*}
\hat{\mathbb E}\Big\{\sup_{\theta\in\Theta}\Big|T^{-1}\int_0^T\tilde\sigma_t|\phi(\theta,X_s)|^{2}ds-T^{-1}g_n^2(\theta)\Big|^{2q}\Big\}\leq L_q(T/n)^q\leq L_q\Delta^{q/2}n^{-q/2}.
\end{align*}
Moreover, for $q>2$, any $\varepsilon>0$, by Markov inequality on sublinear expectation, we get
\begin{align*}
\sum\limits_n{\mathbb V}\Big\{\sup_{\theta\in\Theta}\Big|T^{-1}\int_0^T\tilde\sigma_t|\phi(\theta,X_s)|^{2}ds-T^{-1}g_n^2(\theta)\Big|^{2q}>\varepsilon\Big\}\leq L_q\Delta^{q/2}\sum\limits_n n^{-q/2}/\varepsilon<\infty.
\end{align*}
Thus the Borel-Cantelli lemma of Lemma \ref{BC} yields the required result. Hence, the proof is complete. $\Box$

\begin{lemma}\label{LL} Assume that the random functions $L_n$ satisfy three following conditions:\\
{\rm (C1).} The following equalities hold\\
\begin{align*}
{\mathcal V}\Big\{\underline{L}(\theta)\leq\lim\inf\limits_{n\rightarrow\infty}L_n(\theta)\leq\lim\sup\limits_{n\rightarrow\infty}L_n(\theta)\leq\bar{L}(\theta)\Big\}=1,\\
{\mathbb V}\Big\{\lim\inf\limits_{n\rightarrow\infty}L_n(\theta)=\underline{L}(\theta)\Big\}=1, \quad{\mathbb V}\Big\{\lim\sup\limits_{n\rightarrow\infty}L_n(\theta)=\bar{L}(\theta)\Big\}=1
\end{align*}
uniformly in $\theta\in \Theta$ as $n\rightarrow\infty$.\\
{\rm (C2).} The limiting non-random functions $\underline{L}(\cdot)$ and $\bar{L}(\cdot)$ satisfy
$$\underline{L}(\theta^0)\leq \underline{L}(\theta),\quad\bar{L}(\theta^0)\leq\bar{L}(\theta)\quad\mbox{for all }~\theta\in{\Theta}.$$\\
{\rm (C3).} $\underline{L}(\theta)= \underline{L}(\theta^0)~\mbox{and}~\bar{L}(\theta)=\bar{L}(\theta^0)$ iff $\theta=\theta^0.$

\noindent Then $\hat\theta_n\rightarrow\theta^0$ q.s. as $n\rightarrow\infty$, where

$$L_n(\hat\theta_n)=\min\limits_{\theta\in\Theta}L_n(\theta).$$
\end{lemma}
{\textbf{ Proof}}  From Kasonga \cite[Theorem 2]{Ka}, the claim is easily proved. $\Box$

\section{An example}

In this section, we give an example on the  Ornstein-Uhlenbeck process under sublinear expectations as follows, for $t\in[0,T]$,
\begin{align}\label{Exam}
dX(t)=-\theta^0X(t)dt+dB(t), X(0)=\eta
\end{align}
where $\eta\sim {\cal N}(0,[0.4/\theta^0,0.6/\theta^0]), B(t)\sim {\cal N}(0,[0.5,1]t)$. The solution of the Langevin equation (\ref{Exam}) follows
 $$X^\eta(t):=e^{-\theta^0 t}\Big(\eta+\int_0^te^{\theta^0 u}d B(u)\Big).$$
 By the method of the proof in Lemma \ref{ExI}, together to the idea from Mao \cite[Example 5.1 on page 101]{Msde}, we can deduce that $X^\eta(t)$ is the stationary ergodic process on sublinear expectation space.
In fact, the classical counterpart of (\ref{Exam}) is referred to Bishwal \cite[Example 3.5 (a) on page 95]{Bi}, Cheridito et al. \cite{CKM}, Hu and Nualart \cite{HN}.
Let the O-U process to be observed at equidistant discrete times in equation (\ref{Exam}). Then we have that
\begin{align}\label{Dis}
X(t_{i})=&X(t_{i-1})-\theta^0X(t_{i-1})\Delta_i t+(B(t_{i})-B(t_{i-1}))
\end{align}
with $X(0)=\eta, \Delta_it=T/n$.
Thus, the least squares estimator  of the parameter $\theta^0$ by formula (\ref{LSEa}) is given
\begin{align*}%\label{Exam1}
\hat{\theta}_{n,T}=- \frac{\sum\limits_{i=1}^nX(t_{i-1})(X(t_i)-X(t_{i-1}))/\sigma^2_{t_{i-1}}}{\sum\limits_{i=1}^nX^2(t_{i-1})\Delta_it/\sigma^2_{t_{i-1}}}.
\end{align*}
Let us first introduce the simulation algorithm for $G$-Brownian motion $(B(t), t\in[0,T])$. We consider a random variable $\xi=B({t_i})-B({t_{i-1}})\sim {\cal N}(0,[\underline\sigma^2,\bar\sigma^2]\Delta), i=1,\cdots,n$, we construct an experiment as follows. We take equal-step points $\sigma_k (k=1,\cdots,m)$ such that $\underline\sigma=\sigma_1<\sigma_k<\cdots<\sigma_m=\bar\sigma$. For the $k$th-round random sampling ($k=1,\cdots,m$), $\xi^{kj}_i(i=1,\cdots,n;j=1,\cdots,J)$ are from the classical normal distribution ${\cal N}(0, \sigma_k^2\Delta)$. From (\ref{Dis}), we define $X^{kj}(t_i)$ by
\begin{align}\label{Dis1}
X^{kj}(t_{i})=&X^{kj}(t_{i-1})-\theta^0 X^{kj}(t_{i-1})\Delta +\xi^{kj}_i, ~X_0=\zeta^{kj},
\end{align}
where, $\zeta^{kj}$ are from the $G$-normal distribution ${\cal N}(0,[0.4,0.6])$, $i=1,\cdots,n; j=1,\cdots,J; k=1,\cdots,m$.

Now define
\begin{align}\label{kj}
\hat{\theta}_{n,T}^{kj}=- \frac{\sum\limits_{i=1}^nX^{kj}(t_{i-1})(X^{kj}(t_i)-X^{kj}(t_{i-1}))}{\sum\limits_{i=1}^n(X^{kj}(t_{i-1}))^2\Delta_it}.
\end{align}

   We can obtain the estimator $\bar{\hat\theta}_{n,T}(m,J)$ of $\hat{\mathbb E}[\hat\theta_{n,T}]$ equals to $\max\limits_{1\leq k\leq m}\{\frac{1}{J}\sum\limits_{j=1}^J|\hat\theta_{n,T}^{kj}|\}$, and the estimator $\underline{\hat\theta}_{n,T}(m,J)$ of ${\mathcal E}[\hat\theta_{n,T}^{kj}]$ to $\min\limits_{1\leq k\leq m}\{\frac{1}{J}\sum\limits_{j=1}^J|\hat\theta_{n,T}^{kj}|\}, i=1,\cdots,n$, respectively.

%%%%%%%%%%%%%%%%%%%%%%

\begin{table}
 \textbf{\caption{dependence of LSE on $n$\label{tab1} }}
 \label{table1}
 \begin{tabular}{cccccc}
  \toprule
  $n $ & $1 \times 10^4$  & $ 2 \times 10^4$ &  $3 \times 10^4$   & $4 \times 10^4$  & $5 \times 10^4$     \\
  \midrule
$\bar{\hat{\theta}}_{n,T}(10,2^9)           $ & 1.0313  & 1.0193 &  1.0144  & 1.0089  & 1.0057   \\
$ \underline{\hat{\theta}}_{n,T}(10,2^9)    $ & 1.0104  & 1.0027 &  1.0014  & 0.9972  & 0.9978  \\
$\bar{\hat{\theta}}_{n,T}(10,2^9)- \underline{\hat{\theta}}_{n,T}(10,2^9) $
                                              & 0.0209  & 0.0166 &  0.0130  & 0.0117  & 0.0079  \\
  \bottomrule
 \end{tabular}\\
 Fix $\eta \sim {\cal N}(0,[0.4,0.6]), \Delta t=0.01, m = 10, J = 2^9, T=n\times \Delta t$
 \end{table}
%$\bar{\hat{\theta}}_{n,T}(m,J)$
%$ \underline{\hat{\theta}}_{n,T}(m,J) $

\begin{table}
\textbf{ \caption{dependence of LSE on $J$\label{tab2}}}
 \label{table2}
 \begin{tabular}{cccccc}
  \toprule
  $J$ & $2^{3}$  & $2^4$  & $2^5$ &  $2^6$   & $2^7$      \\
  \midrule
$\bar{\hat{\theta}}_{5 \times 10^3,T}(10,J)$           & 1.1108 & 1.0955 &  1.0903   & 1.0780  & 1.0672     \\
$ \underline{\hat{\theta}}_{5 \times 10^3,T}(10,J) $   & 0.9718 & 0.9713 &  0.9879   & 0.9888  & 0.9995     \\
 $\bar{\hat{\theta}}_{5 \times 10^3,T}(10,J)-\underline{\hat{\theta}}_{5 \times 10^3,T}(10,J)$
                                              & 0.1390 & 0.1241 &  0.1023   & 0.0892  & 0.0677    \\
  \bottomrule
 \end{tabular}\\
 Fix $\eta \sim {\sl N}(0,[0.4,0.6]), T = 50, m = 10, \Delta t = 0.01$
\end{table}

Let $\theta^0=1,T=0.01 n$ in (\ref{Exam}). Note that $\underline{\sigma}^2=0.5,\bar\sigma^2=1$. If we fix $m=10,J=2^9$, for different values $n$, by (\ref{Dis1})-(\ref{kj}), we obtain the numerical consequence of lower and upper expectation pair  $(\bar{\hat\theta}_{n,T}(m,J), \underline{\hat\theta}_{n,T}(m,J)$ of the parameter LSE $\hat\theta_{n,T}$  as  Table 1.
 If we fix $m=10,n=10^3, T=50$, for different values $J$, by (\ref{Dis1})-(\ref{kj}), we obtain the numerical consequence of lower and upper expectation pair  $(\bar{\hat\theta}_{n,T}(m,J), \underline{\hat\theta}_{n,T}(m,J)$ of the parameter LSE $\hat\theta_{n,T}$  as  Table 2.

 Table \ref{tab1} and \ref{tab2} show that the LSE of the unknown parameter $\theta^0$ strongly converge to the true value quasi-surely and the errors of convergence change more and more small as $n\rightarrow\infty$, $J\rightarrow\infty$, respectively. Thus, our main theorem is verified numerically.

\section{Conclusion}

  In reality, we often are faced with probability and Knightian uncertainties, respectively. A real system with ambiguity can often be characterized by $G$-Brownian motion by using Peng's sublinear expectation framework. Thus we need to consider a stochastic differential equation driven by $G$-Brownian motion. Moreover, if a $G$-SDE has unknown parameters then we are asked to estimate it. In this paper, we gives an estimator of LSE of the unknown parameter, and proved the strong convergence theorem. Moreover, to complete the proof of theorem, we prepare several lemmas. Finally, we discuss an example of the Ornstein-Uhlenbeck process, where maxmum and minimum estimators are given. By numerical simulation, we show the behaviour of the LSE of the parameter of $G$-SDE, which show the strong convergence numerically.

     Since an observer has multi-priors, the least squares estimator of the parameter often take multi-values. We know that a Knightian uncertainty averse economic agent looks for a robust portfolio make-decision in finance. If we assume that the observer of the stochastic system with distribution uncertainty is averse to uncertainty of the estimation error. This will cause a problem of how to determine the only estimation value of the unknown parameter, which depends the observer's attitude to multi-estimation values. However, in this paper, we only study the LSE in (\ref{LSEa}) without embodying the observer's attitude to the multi-estimation values. In future, we shall study the strong consistency and the limit distribution of {\sl robust least squares estimator} of the unknown parameter of $G$-SDE.

\vskip12pt

{\bf Acknowledgements}  This paper is accepted for publication in Acta Mathematica Scientia (Chinese Series).

\end{document}